\documentclass[11pt]{amsart}

\usepackage{amssymb,bbm,amscd}
\usepackage{graphicx}
\usepackage{a4wide}

\theoremstyle{plain}

\theoremstyle{definition}

\newcommand{\ts}{\hspace{0.5pt}}
\newcommand{\nts}{\hspace{-0.5pt}}

\newcommand{\ZZ}{\mathbb{Z}}
\newcommand{\RR}{\mathbb{R}\ts}

\newcommand{\NN}{\mathbb{N}}

\begin{document}

\title[Factors of substitution systems]{Examples of substitution 
systems and their factors\vspace*{5ex}}

\author{Michael Baake}
\author{Franz G\"{a}hler}

\address{Fakult\"{a}t f\"{u}r Mathematik, Universit\"{a}t Bielefeld,\newline
\hspace*{\parindent}Postfach 100131, 33501 Bielefeld, Germany}
\email{$\{$mbaake,gaehler$\}$@math.uni-bielefeld.de}

\author{Uwe Grimm}
\address{Department of Mathematics and Statistics,
The Open University,\newline 
\hspace*{\parindent}Walton Hall, Milton Keynes MK7 6AA, United Kingdom}
\email{uwe.grimm@open.ac.uk}

\begin{abstract}
  The theory of substitution sequences and their higher-dimensional
  analogues is intimately connected with symbolic dynamics. By
  systematically studying the factors (in the sense of dynamical
  systems theory) of a substitution dynamical system, one can reach a
  better understanding of spectral and topological properties. We
  illustrate this point of view by means of some characteristic
  examples, including a rather universal substitution in one dimension
  as well as the squiral and the table tilings of the plane.
\end{abstract}

\maketitle
\thispagestyle{empty}

\bigskip
\bigskip

\section{Introduction}

Closed subshifts over finite alphabets are much studied dynamical
systems \cite{LMBook,Kit}. Their understanding in one dimension is a
cornerstone of the theory of (symbolic) dynamical systems, while
systems in higher dimensions show many new phenomena
\cite{Schmidt}. It is probably fair to say that in the latter case,
despite great effort, the open questions still prevail.  The class of
symbolic dynamical systems (over finite alphabets) is special also in
the sense that they possess finitely many (non-periodic) factors (up
to isomorphism), where the term `factor' refers to symbolic dynamical
systems that are the image of a homomorphism which commutes with the
shift action (and not to finite subwords).

Nevertheless, the factors are rarely used explicitly to unravel
(details of) the structure of the dynamical system. For certain
aspects, however, they carry relevant information, and can be employed
both for structural insight and for concrete calculations.  This is
particularly true for the spectral theory of symbolic dynamics. The
latter comes in two flavors. The traditional point of view is via the
\emph{dynamical} spectrum \cite{Q}, which analyzes the induced shift
action on the Hilbert space of square integrable functions on the
shift space (relative to an invariant probability measure). An
alternative approach works with the \emph{diffraction} spectrum
\cite{Hof,M,BGrev} of a typical sequence in the space, which was
originally motivated by the physical process of kinematic diffraction
\cite{Cow}. Its relevance increased with the discovery of
quasicrystals by Shechtman et al.~\cite{Danny}, particularly for
systems with pure point (or Bragg) diffraction. Recent evidence
\cite{Withers} suggests that also continuous spectral components are
practically relevant.

It is well understood by now that (and how) the two types of spectra
are equivalent in the pure point case \cite{LMS,M,BL}, while it has
long been known that the dynamical spectrum is generally richer in the
presence of continuous components \cite{EM,BE}. Nevertheless, in many
of the classic examples, the complete dynamical spectrum can be
reconstructed from the diffraction spectrum of the system and of some
of its factors. This can be made more precise for systems of symbolic
dynamics as well as for Delone dynamical systems with finite local
complexity \cite{BEL}. We will see examples of this type below,
selected from the class of substitution systems with integer inflation
factor; see \cite{Fre02,Nat1,Nat2,LM,FS} for background material.

Yet another, and rather recent, aspect is the topological structure of
tiling spaces; see \cite{Sadun} and references therein. Here, the
substitution or the inflation structure can effectively be used to
calculate the \v{C}ech cohomology, for instance via the methods
introduced in \cite{AP}. This provides concrete results, but their
(geometric) interpretation is still difficult. This is another
instance where the careful inspection of the factors can help to
understand the detailed structure of the cohomology groups. We will
explain this in some detail for a recently analyzed example, the
so-called squiral tiling of the plane \cite{GS,BG12}.

Our approach in this article is example-oriented, wherefore we opted
for a somewhat informal presentation, in the interest of better
readability. We assume the reader to be familiar with the basic
concepts of (symbolic) dynamics, compare \cite{LMBook,Rob04}, and with
some results from their spectral theory, see \cite{Q} for more.  We
begin with an example in one dimension that is deliberately designed
to have interesting factors, and a rather rich spectrum as a result of
this. We continue with some remarks on substitution factors with
maximal pure point spectrum, before we embark on a more detailed
analysis of the squiral tiling \cite{GS,BG12} and its substitution
factors.  This is followed by a similar (though still somewhat
preliminary) analysis of the classic table tiling \cite{Rob99}.

\section{A fun example: The `universal' morphism}\label{sec:fun}

Undeniably, the best studied substitution rule (on the binary alphabet
$\{a,b\}$) is given by $\varrho^{}_{\mathrm{F}}\!:\, a\mapsto ab,\,
b\mapsto a$ and known as the Fibonacci substitution (or morphism). It
defines a unique hull $\mathbb{X}_{\mathrm{F}} \subset
\{a,b\}^{\mathbb{Z}}$ that is strictly ergodic as a dynamical system
under the action of the (two-sided) shift; see \cite{AS,Q,PF,TAO} and
references therein. It has pure point dynamical as well as pure point
diffraction spectrum. In a similar way, one can treat all (repetitive)
Sturmian sequences, a common feature being that they can all be
described as model sets (also known as cut and project sets); see
\cite{Moody00} as well as \cite{PF} for a survey. The equivalence of
the two types of spectra is well understood for systems of this kind
\cite{LMS,BL,BLM}.

The relation between the dynamical and the diffraction spectrum is
more complicated in the presence of continuous spectral components.
This was first pointed out in \cite{EM} for the example of the
Thue-Morse substitution, where the dynamical spectrum is richer. In
particular, a substantial part of the pure point spectrum does not
show up in the diffraction measure of the Thue-Morse chain, while it
can be recovered from the diffraction of the period doubling
chain. The hull of the latter is a factor of the Thue-Morse system
\cite{Sadun,BGG}. An analogous phenomenon was described in \cite{BE}
for a system with close-packed dimers on $\mathbb{Z}$, and is known to
also occur for the Rudin-Shapiro chain \cite{TAO}. The purpose of this
section is to illustrate these connections with a single substitution
rule that entails a mixed spectrum and possesses most of the above
examples as factors.

Let us consider the alphabet $\{a,b,c,d,\bar{a},\bar{b},\bar{c},\bar{d}\}$
and define the primitive substitution
\begin{equation}\label{eq:unimorph}
  \varrho : \quad
  \begin{array}{c@{\;}c@{\;}ccc@{\;}c@{\;}c}
    a & \mapsto & a \bar{b} & \quad & \bar{a} & \mapsto & \bar{a} b \\
    b & \mapsto & a \bar{d} & \quad & \bar{b} & \mapsto & \bar{a} d \\
    c & \mapsto & c\ts \bar{d} & \quad & \bar{c} & \mapsto & \bar{c}\ts  d \\
    d & \mapsto & c\ts \bar{b} & \quad & \bar{d} & \mapsto & \bar{c}\ts  b 
  \end{array}
\end{equation}
which results in a unique hull, for instance via the two-sided sequence
\[
  w = \ldots a\bar{b}\bar{c}b\bar{c}da\bar{d}\bar{c}dc\bar{b}a\bar{b}\bar{c}b
   |a\bar{b}\bar{a}d\bar{a}bc\bar{b}\bar{a}ba\bar{d}c\bar{d}\bar{a}d \ldots
\]
which is a fixed point under $\varrho^{2}$ with legal seed $b|a$. The
corresponding hull $\mathbb{X}$ is obtained as the closure of the
two-sided shift orbit of $w$ in the product topology. Before
we comment on the spectral structure, let us identify some relevant
factors.

Upon identifying letters with their barred copies, one obtains
the \emph{quaternary} Rudin-Shapiro substitution \cite{AS,Q}
\[
    \varrho^{}_{\mathrm{RS}} : \quad
    a \mapsto ab \, , \quad b \mapsto ad \, , \quad
    c \mapsto cd \, , \quad d \mapsto cb\, ,
\]
with induced fixed point $w^{}_{\mathrm{RS}}$ of
$\varrho^{2}_{\mathrm{RS}}$, with legal seed $b|a$.  This leads to the
classic \emph{binary} Rudin-Shapiro chain via the map $\varphi$
defined by $\varphi (a) = \varphi (b) = 1$ and $\varphi (c) = \varphi
(d) = \bar{1}$. For complexity results, we refer to \cite{AS93}. The
corresponding hulls fail to be palindromic \cite{Apal,Bpal}, which has
interesting consequences for the spectral theory of associated
Schr\"{o}dinger operators. For further results on palindromic systems,
we refer to \cite{ABCD} and references therein.

The binary and the quaternary Rudin-Shapiro hulls are \emph{mutually
  locally derivable} (MLD); see \cite{B02} and references therein for
the concept.  In the symbolic setting, MLD just means that there are
sliding block maps with local support \cite{LM} in both directions of
the derivation.  The non-trivial direction of this claim follows from
the observation that $\bar{1}\bar{1}\bar{1}\bar{1}$ is the longest
$1$-free subword of the binary image and has the unique preimage
$dcdc$ under $\varphi$. This determines the two cosets of $\ZZ$ modulo
$2$ that are occupied by the letters $a,c$ and $b,d$,
respectively. Due to repetitivity, the subword
$\bar{1}\bar{1}\bar{1}\bar{1}$ occurs with bounded gaps, wherefore
this defines a \emph{local} derivation rule. Consequently, both define
dynamical systems with the same spectrum, which is known to comprise
the pure point part $\ZZ[\frac{1}{2}]$ and an absolutely continuous
component, the latter being Lebesgue measure with multiplicity $2$;
see \cite{Q} for details.

Let us go one step back, and define a sliding block map 
on $\mathbb{X}_{\mathrm{RS}}$ that is induced by
\[
\begin{split}
    \chi(ab) = \chi(cd) = A\, , \quad 
    \chi(ad) = \chi(cb) = B\, ,\\
    \chi(ba) = \chi(dc) = C\, ,\quad
    \chi(da) = \chi(bc) = D\, .
\end{split}
\]
Its action on $\mathbb{X}_{\mathrm{RS}}$ is given by $u\mapsto
u^{}_{\chi}$ with $u^{}_{\chi}(i):=\chi\bigl(u(i)u(i+1)\bigr)$.  This
way, one induces the primitive substitution
\[
   \varrho^{}_{\mathrm{T}}:\quad
   A\mapsto AC\, ,\quad
   B\mapsto AD\, ,\quad
   C\mapsto BD\, ,\quad
   D\mapsto BC\, ,
\]
which shows a Toeplitz structure (see below) that is somewhat similar
to that of the paper folding substitution; compare \cite{AB,AMF,TAO}. 

Let $w^{}_{\mathrm{T}}$ be the fixed point of
$\varrho^{2}_{\mathrm{T}}$ with legal seed $C|A$, which is the image
of $w^{}_{\mathrm{RS}}$ under the sliding block map induced by $\chi$.
This bi-infinite word leads to a partition of $\mathbb{Z}$ into four
subsets $\varLambda_{\alpha}=\{i\in\mathbb{Z}\mid
w^{}_{\mathrm{T}}(i)=\alpha\}$ for $\alpha\in\{A,B,C,D\}$. The fixed
point property then results in the equations
\[
   \begin{split}
   \varLambda^{}_{A}&=4\mathbb{Z}\, , \qquad 
   \varLambda^{}_{B}=4\mathbb{Z}+2\, , \\
   \varLambda^{}_{C}&=(8\mathbb{Z}+1)\cup (16\mathbb{Z}+11)\cup 
   (4\varLambda^{}_{C}+3)\, ,\\
   \varLambda^{}_{D}&=(8\mathbb{Z}+5)\cup (16\mathbb{Z}+3)\cup 
   (4\varLambda^{}_{D}+3)\, ,
   \end{split}
\]
where $\varLambda^{}_{A}\cup\varLambda^{}_{B}=2\mathbb{Z}$ and
$\varLambda^{}_{C}\cup\varLambda^{}_{D}=2\mathbb{Z}+1$ was used to
simplify the relations. Due to the decoupling, one can calculate the
solutions via iteration, which gives
\[
    \begin{split}
    \varLambda^{}_{C} & = \{-1\} \cup
                         \bigcup_{n\ge 0}\bigl(2\cdot 4^{n+1}\mathbb{Z}
                         + (2\cdot 4^{n}-1) \bigr) \cup
                         \bigcup_{n\ge 1}\bigl(4^{n+1}\mathbb{Z}
                         + (3\cdot 4^{n}-1) \bigr), \\
    \varLambda^{}_{D} & = \bigcup_{n\ge 0}\bigl(2\cdot 4^{n+1}\mathbb{Z}
                         + (6\cdot 4^{n}-1) \bigr) \cup
                         \bigcup_{n\ge 1}\bigl(4^{n+1}\mathbb{Z}
                         + (4^{n}-1) \bigr).
    \end{split}
\]
Note that the singleton set $\{-1\}$ has to be added to
$\varLambda^{}_{C}$ because it is not contained in any of the lattice
cosets (although it is in the $2$-adic closure of either of the four
unions).  The only other fixed point of $\varrho^{2}_{\mathrm{T}}$
differs from $w^{}_{\mathrm{T}}$ precisely in
$w^{}_{\mathrm{T}}(-1)$. It can thus be described by moving the
singleton set $\{-1\}$ from $\varLambda^{}_{C}$ to
$\varLambda^{}_{D}$. Our explicit coordinatization is a result of a
coincidence in the sense of Dekking \cite{D}. It implies that
$\varrho^{}_{\mathrm{T}}$ defines a dynamical system with pure point
dynamical spectrum, which is $\mathbb{Z}[\frac{1}{2}]$. The latter
coincides with the pure point part of the spectrum of the
Rudin-Shapiro system. This also follows from standard results; compare
\cite{Nat1} and references therein.

Returning to the original substitution $\varrho$ of
Eq.~\eqref{eq:unimorph} and mapping all ordinary letters to $1$ and
all barred ones to $\bar{1} = -1$, one induces the Thue-Morse (or
Pruhet-Thue-Morse) substitution
\[
    \varrho^{}_{\mathrm{TM}} : \quad 1 \mapsto 1 \bar{1} \, , \quad
    \bar{1} \mapsto \bar{1} 1 \, ,
\]
formulated on the alphabet $\{ 1, \bar{1} \}$; compare \cite{AS} for
background. Our original fixed point $w$ is mapped to
$w^{}_{\mathrm{TM}}=\varrho^{2}_{\mathrm{TM}}(w^{}_{\mathrm{TM}})$
with legal seed $1|1$.  The dynamical spectrum of the corresponding
hull $\mathbb{X}_{\mathrm{TM}}$ under the shift action comprises a
pure point part (namely $\mathbb{Z}[\frac{1}{2}]$) and a singular
continuous one, the latter leading to an explicit representation as a
Riesz product \cite{Kaku,Q,BGG}.

Defining the sliding block map $\psi$ on $\mathbb{X}_{\mathrm{TM}}$
via $\psi(w) (i) = - w(i) w(i+1)$, one induces another classic
substitution,
\[
    \varrho^{}_{\mathrm{pd}} : \quad
    1 \mapsto 1\bar{1} \, , \quad \bar{1} \mapsto 11 \, ,
\]
which is known as the period doubling substitution; compare
\cite{AS}. The image of $w^{}_{\mathrm{TM}}$ under $\psi$ is
$w^{}_{\mathrm{pd}}=\varrho^{2}_{\mathrm{pd}}(w^{}_{\mathrm{pd}})$,
with legal seed $\bar{1}|1$.  A coordinatization of
$w^{}_{\mathrm{pd}}$ via a partition of $\mathbb{Z}$ can be done in a
similar way as discussed above for $\varrho^{}_{T}$; compare
\cite{BM}.  Via the corresponding diffraction measure, this provides an
alternative way to show that $\mathbb{Z}[\frac{1}{2}]$ is the
dynamical spectrum of the associated dynamical system. It exhausts the
pure point part of the dynamical spectrum of the Thue-Morse system.

Having analyzed these (selected) factors of our original substitution
$\varrho$, we can conclude that the latter has mixed spectrum with all
three spectral types being present. The maximal equicontinuous (or
Kronecker) factor is the dyadic solenoid $\mathbb{S}^{1}_{2}$, which
also emerges via the `torus parametrization' of the period doubling
chain or the maximal model set factor of the Rudin-Shapiro chain.
Let us explain this type of connection in a little more detail.

\section{Substitution factors with maximal 
   pure point spectrum}\label{sec:Kron}

 The Thue-Morse dynamical system has a dynamical spectrum of mixed
 type, with pure point and singular continuous components.  As was
 pointed out in \cite{EM}, the diffraction spectrum only detects part
 of it. In fact, it only shows the trivial part of the pure point
 spectrum, $\ZZ$, while the rest is missing (or hidden). However, the
 period doubling system is a factor that is pure point, with dynamical
 spectrum $\ZZ[\frac{1}{2}]$, so that the diffraction of Thue-Morse
 together with that of period doubling covers the entire dynamical
 spectrum of Thue-Morse.

 In fact, the period doubling system is a factor that lies between
 Thue-Morse and its Kronecker factor, the latter being the dyadic
 solenoid $\mathbb{S}^{1}_{2}$ \cite{Nat1,Sadun}. This happens to be
 the `torus' for the period doubling system in the language of model
 sets \cite{BLM}, and the mapping from period doubling to the solenoid
 is 1-to-1 almost everywhere \cite{BM,BLM}.

 The system of close-packed dimers on the line \cite{BE} provides an
 example where the dynamical spectrum has pure point and absolutely
 continuous components, as does Rudin-Shapiro. In both cases, the
 dynamical spectrum is richer, but the `missing' parts again can be
 recovered from the analysis of a single factor.  A similar situation
 shows up in many other examples, though not in all. Still, it is
 worth looking at this relation in some more detail.

A natural generalization of the Thue-Morse substitution, in the spirit
of \cite{Kea}, is
\[
    \varrho^{(k,\ell)}_{\ts\mathrm{gTM}} : \;
    \begin{array}{r@{\;}c@{\;}l}
    a & \mapsto & a^{k}b^{\ell} \\ b & \mapsto & b^{k}a^{\ell}
    \end{array}
\]
with $k,\ell\in\NN$, where $k=\ell=1$ is the classic TM substitution.
The substitution matrix reads $\left(\nts\begin{smallmatrix} k & \ell
    \\ \ell & k\end{smallmatrix}\nts\right)$, with eigenvalues $k\pm
\ell$. This two-parameter family of constant length substitutions
shares many properties with its classic ancestor; see \cite{BGG} for
a detailed analysis.  In particular, there is a sliding block map that
works for the entire family, and induces a maximal model set factor of
substitution type.  The latter is a generalization of the period
doubling substitution, namely
\[
    \varrho^{(k,\ell)}_{\ts\mathrm{gpd}} : \;
    \begin{array}{r@{\;}c@{\;}l}
    a & \mapsto & ub \\ b & \mapsto & ua
    \end{array}
\]
with the finite word $u=b^{k-1}ab^{\ell-1}$.  This is a constant length
substitution with a coincidence in the sense of Dekking \cite{D}, and
hence a model set \cite{LM}. This factor then has a torus
parametrization \cite{BLM} which is an almost everywhere 1-to-1 map
onto a one-dimensional dyadic solenoid, which we denote by
$\mathbb{S}^{1}_{2}$ as before.

Let us now discuss some two-dimensional examples, which are
considerably more complex.

\section{The squiral tiling and its factors}\label{sec:squiral}

\begin{figure}
  \centerline{\includegraphics[width=0.9\textwidth]{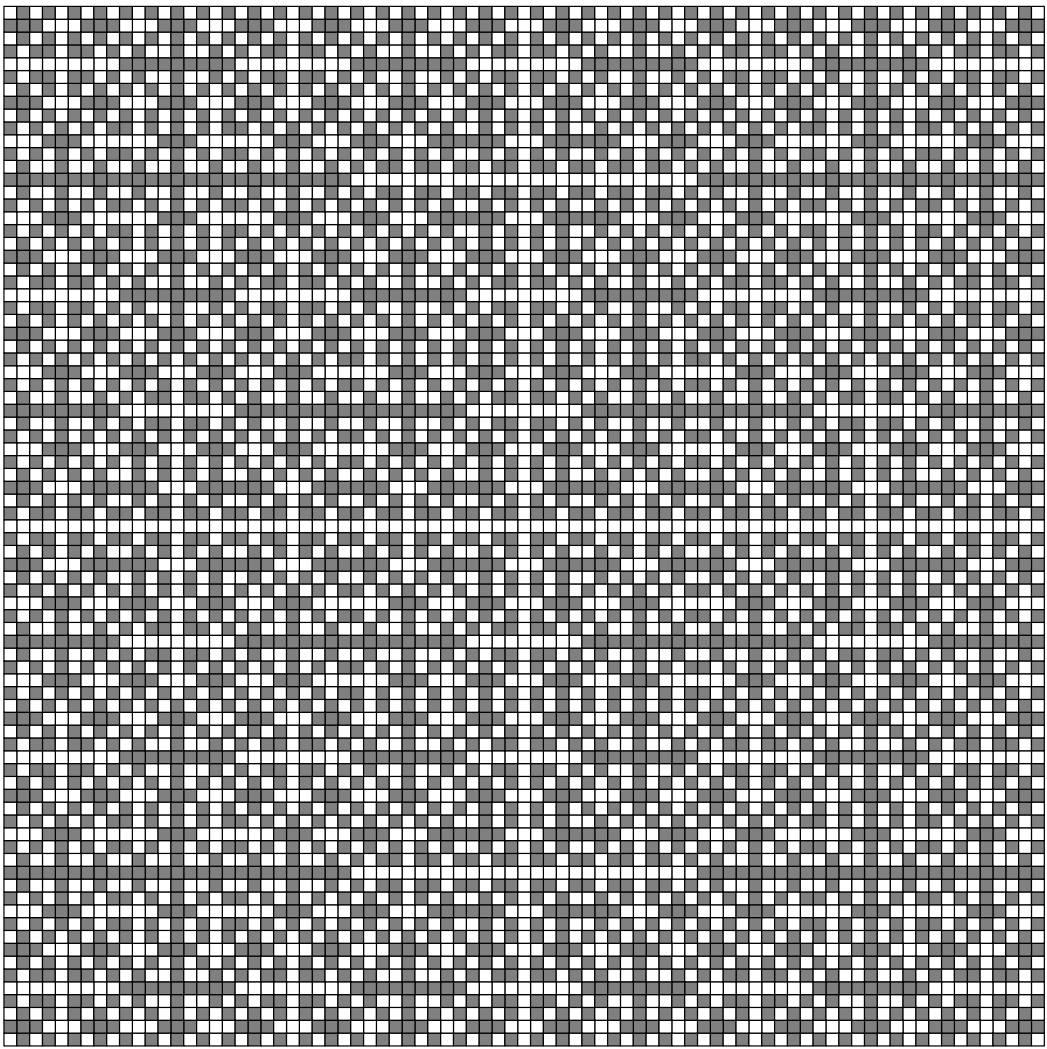}}
\caption{\label{fig:blocktil} A square-shaped patch of the squiral
 tiling with full $D_{4}$-symmetry, as obtained from the 
 block substitution \eqref{eq:blocksub}. For this figure, the
 reference point of the prototiles is placed in their centers,
 so that the patch is the central part of a fixed point under the
 corresponding block substitution.}
\end{figure}

The \emph{squiral tiling} appears in \cite[Fig.~10.1.4]{GS}, where it
was constructed as a simple example for a tiling of the plane by one
prototile (and its mirror image) with infinitely many edges. It is
obtained from the primitive inflation rule
\begin{equation}\label{eq:sqinfl}
   \raisebox{-8ex}{\includegraphics[width=0.4\textwidth]{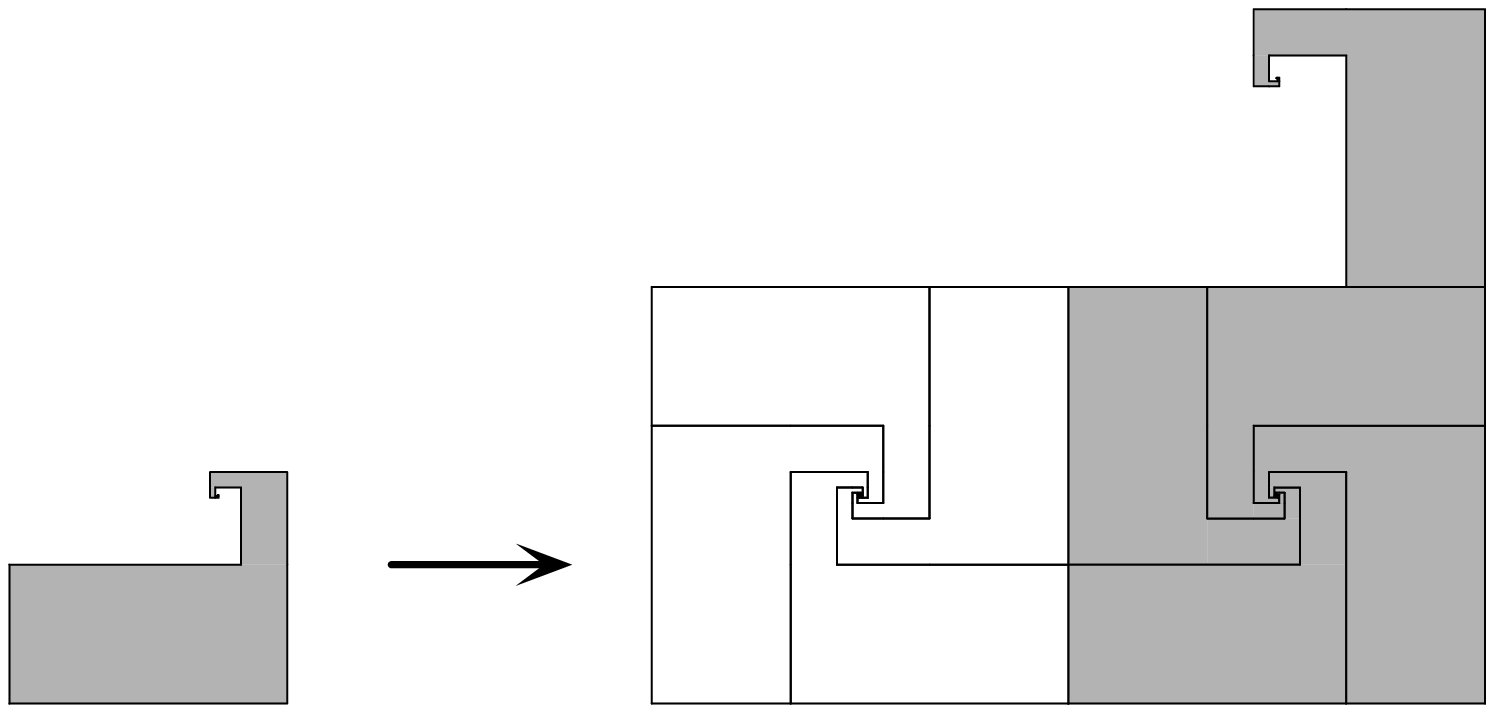}}
\end{equation}
which has the (linear) integer inflation factor $3$. When one
distinguishes the two chiralities, it is obvious that the (colored)
tiling is equivalent (in the sense of mutual local derivability
\cite{BSJ,B02}) to a coloring of the square lattice that emerges from
the block substitution
\begin{equation}\label{eq:blocksub}
  \raisebox{-6ex}{\includegraphics[width = 0.6 \textwidth]{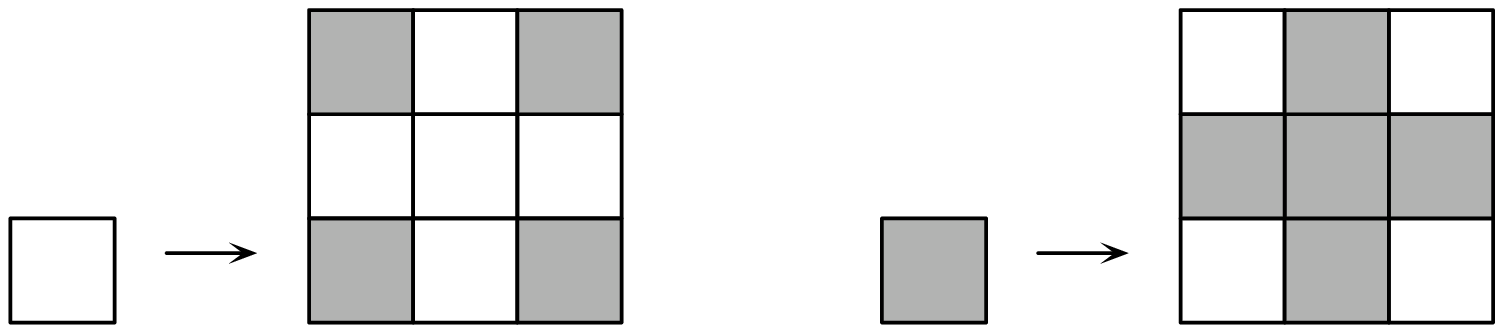}}
\end{equation}
This block substitution is bijective and of constant length in the
terminology of \cite{Nat1}. In this formulation, it was analyzed in
detail in \cite{BG12}. A larger patch is shown in Figure~\ref{fig:blocktil}.

It follows from standard arguments \cite{Nat1,Nat2} that the dynamical
spectrum of the squiral contains
$\ZZ[\frac{1}{3}]\times\ZZ[\frac{1}{3}]$ as its pure point
part. Moreover, the analysis of \cite{BG12} shows that the remaining
part of the spectrum is purely singular continuous. The corresponding
diffraction measure is explicitly known, and can be represented as a
two-dimensional Riesz product.

As the tiling is bijective \cite{Nat1}, the block substitution
(\ref{eq:blocksub}) commutes with the exchange of the two colors.
This implies that the hull of the squiral is also invariant under the
exchange of the colors: with every tiling in the hull, also the
tiling with inverted colors is in the hull. This is the deeper reason
why the projection onto the underlying two-dimensional $3$-adic
solenoid $\mathbb{S}^{2}_{3}$ is at least 2-to-1, so that the
(dynamical or diffraction) spectrum cannot be pure point. In
particular, the squiral tiling cannot originate from a model set.

Our next step therefore is to find a factor in which pairs of tilings
with inverted colors are identified. This is achieved with a sliding
block map, which maps $2\times2$-blocks of tiles to some new symbols.
The squiral tiling admits 14 distinct legal $2\times2$-patterns, as
only the two blocks with four equal symbols are not allowed. In the
following block map (in which we write the two colors as $1$ and
$\bar{1}=-1$), pairs of blocks with inverted colors are identified,
but nothing else,
\begin{equation}
\begin{gathered}
\pm\left[\begin{matrix} 1 & 1 \\ \bar{1} & 1 \end{matrix}\right]\mapsto a,\quad
\pm\left[\begin{matrix} 1 & 1 \\ 1 & \bar{1} \end{matrix}\right]\mapsto b,\quad
\pm\left[\begin{matrix} 1 & \bar{1} \\ 1 & 1 \end{matrix}\right]\mapsto c,\quad
\pm\left[\begin{matrix} \bar{1} & 1 \\ 1 & 1 \end{matrix}\right]\mapsto d,\\
\pm\left[\begin{matrix} 1 & 1 \\ \bar{1} & \bar{1} \end{matrix}
\right]\mapsto e,\quad
\pm\left[\begin{matrix} 1 & \bar{1} \\ 1 & \bar{1} \end{matrix}
\right]\mapsto f,\quad
\pm\left[\begin{matrix} 1 & \bar{1} \\ \bar{1} & 1 \end{matrix}
\right]\mapsto g.
\end{gathered}
\label{eq:blockmap}
\end{equation}
It is easy to see that any tiling in the image lifts to exactly two
squiral tilings. There is one position where we can chose whether the
lift has a $1$ or a $\bar{1}$ there; all other tiles in the lift are
then determined.  Our factor map is therefore uniformly 2-to-1 on the
entire hull of all squiral tilings. This is completely analogous to
the hull of the (generalized) Thue-Morse chains from above, which
project uniformly 2-to-1 to the hull of the (generalized)
period doubling chains, under a similar sliding block map \cite{BGG}.

The block map (\ref{eq:blockmap}) induces a primitive substitution on the 
tilings in the factor space,
\begin{equation}
\begin{gathered}
a\mapsto\left[\begin{matrix} g&g&a \\ d&c&g \\ a&b&g \end{matrix}\right],\quad
b\mapsto\left[\begin{matrix} f&f&b \\ d&c&g \\ a&b&g \end{matrix}\right],\quad
c\mapsto\left[\begin{matrix} f&f&c \\ d&c&e \\ a&b&e \end{matrix}\right],\quad
d\mapsto\left[\begin{matrix} g&g&d \\ d&c&e \\ a&b&e \end{matrix}\right],\\
e\mapsto\left[\begin{matrix} g&g&e \\ d&c&e \\ a&b&e \end{matrix}\right],\quad
f\mapsto\left[\begin{matrix} f&f&f \\ d&c&g \\ a&b&g \end{matrix}\right],\quad
g\mapsto\left[\begin{matrix} g&g&g \\ d&c&g \\ a&b&g \end{matrix}\right].
\end{gathered}
\label{eq:factorsubst}
\end{equation}
Obviously, the resulting tilings have a Toeplitz structure; compare
\cite{AB}. They contain lattice-periodic subsets, which implies that
they are model sets \cite{LM}.  As such, they necessarily have pure
point dynamical and diffraction spectra; compare \cite{LMS}.  In turn,
the squiral tilings, whose hull is a 2-to-1 cover of a pure point
factor, cannot have pure point spectrum. In fact, with suitable
scattering strengths ($+1$ and $-1$ on the two symbols), their
diffraction spectrum is purely singular continuous \cite{BG12}.

The tilings generated by the substitution (\ref{eq:factorsubst}) form
the maximal model set factor of the squiral tiling, which we
henceforth call $F_\text{max}$. This factor, in turn, has a
two-dimensional, $3$-adic solenoid $\mathbb{S}_3^2$ as factor, onto
which it projects 1-to-1 almost everywhere via the torus
parametrization of general model sets \cite{BLM}. We shall now analyze
in more detail the set of tilings where this projection fails to be
1-to-1.

For this purpose, we note that the level-$n$ supertiles of the
substitution (\ref{eq:factorsubst}) have the following block structure
\begin{equation}
\begin{gathered}
a\mapsto\left[\begin{matrix} G & a \\ X & G^T \end{matrix}\right],\quad
b\mapsto\left[\begin{matrix} F & b \\ X & G^T \end{matrix}\right],\quad
c\mapsto\left[\begin{matrix} F & c \\ X & E^T \end{matrix}\right],\quad
d\mapsto\left[\begin{matrix} G & d \\ X & E^T \end{matrix}\right],\\
e\mapsto\left[\begin{matrix} G & e \\ X & E^T \end{matrix}\right],\quad
f\mapsto\left[\begin{matrix} F & f \\ X & G^T \end{matrix}\right],\quad
g\mapsto\left[\begin{matrix} G & g \\ X & G^T \end{matrix}\right],
\end{gathered}
\label{eq:supertiles}
\end{equation}
where $X$ is a square block of dimension $(3^n-1)$ (which is the same
for all seven symbols), $F$ and $G$ are rows with entries $f$ and $g$,
while $E^T$ and $G^T$ are columns with entries $e$ and $g$,
respectively.

We can now arrange two supertiles such that the separating line between them 
passes near the origin. Two blocks $X$ are then separated  by a line in which 
all symbols are the same. There are four possibilities, two with a horizontal 
and two with a vertical separation line,
\begin{equation}
\left[\begin{matrix} X \\ F \\ X \end{matrix}\right],\quad
\left[\begin{matrix} X \\ G \\ X \end{matrix}\right],\quad
\left[\begin{matrix} X & E^T & X \end{matrix}\right],\quad
\left[\begin{matrix} X & G^T & X \end{matrix}\right].
\end{equation}
Taking together all ways in which the boundaries of the supertiles can
tend to infinity as the supertile order is increased, we obtain four
one-dimensional model sets, two of which are arranged horizontally,
and two vertically. Projected to $\mathbb{S}^{2}_{3}$, the translation
orbits of the horizontal pair project to a single translation orbit of
a one-dimensional sub-solenoid of type $\mathbb{S}^{1}_{3}$, and so do
the translation orbits of the vertical pair.

If a quartet of infinite order supertiles is placed such that the
common corner remains near the origin, the following configurations
are obtained,
\begin{equation}
\begin{gathered}
\left[\begin{matrix} X&E^T&X \\ G&a&F \\ X&G^T&X \end{matrix}\right],\quad
\left[\begin{matrix} X&E^T&X \\ F&b&G \\ X&G^T&X \end{matrix}\right],\quad
\left[\begin{matrix} X&G^T&X \\ F&c&G \\ X&E^T&X \end{matrix}\right],\quad
\left[\begin{matrix} X&G^T&X \\ G&d&F \\ X&E^T&X \end{matrix}\right],\\
\left[\begin{matrix} X&E^T&X \\ G&e&G \\ X&E^T&X \end{matrix}\right],\quad
\left[\begin{matrix} X&G^T&X \\ F&f&F \\ X&G^T&X \end{matrix}\right],\quad
\left[\begin{matrix} X&G^T&X \\ G&g&G \\ X&G^T&X \end{matrix}\right].
\end{gathered}
\label{eq:specpoints}
\end{equation}
Each of these configurations represents a translation orbit in $F_\text{max}$, 
and all these orbits are projected to a single orbit in 
$\mathbb{S}_3^2$. The configurations shown in (\ref{eq:specpoints}) form
the seeds of the seven fixed points of the substitution 
(\ref{eq:factorsubst}).

This structure of the hull is in line with the Artin-Mazur dynamical
zeta function \cite{Ruelle} of the substitution action on the hull of
$F_\text{max}$. On the one hand, the zeta function is defined in terms
of the periodic points in the hull under the substitution,
\begin{equation}
   \zeta(z) \, = \, \exp\left(\sum_{m=1}^\infty 
   \frac{a_m}{m} z^m\right) \, = \, \prod_{m=1}^{\infty}
   (1-z^{m})_{}^{-c_{m}} ,
\label{eq:zetaperpoints}
\end{equation}
where $a_m$ is the number of points in the hull that are invariant
under an $m$-fold substitution. Likewise, the exponents $c_m$ in the
Euler product are the cycle numbers, which follow from the $a_m$ via
\[
    c_{m} \, = \, \frac{1}{m} \sum_{d|m} \mu(\tfrac{d}{m})\ts a_{d} \ts ,
\]
where $d$ runs through the divisors of $m$ and $\mu$ denotes the
M\"{o}bius function of elementary number theory; see
\cite[Ch.~6.4]{LMBook} for a detailed exposition.  Note that if the
hull consists of several components for which the periodic points can
be counted separately, the total zeta function is obtained as the
product of the partial zeta functions. In our case, according to the
analysis above, $F_\text{max}$ consists of one copy of
$\mathbb{S}^{2}_{3}$, two extra copies of one-dimensional solenoids
$\mathbb{S}^{1}_{3}$ (above those points where the projection to
$\mathbb{S}^{2}_{3}$ is 2-to-1), and four extra fixed points above the
origin of $\mathbb{S}^{2}_{3}$.

This structure of the hull and the dynamical zeta function can be
confirmed by computing the dynamical zeta function by a second method,
which is a by-product of the computation of the \v{C}ech cohomology of
the hull via the Anderson-Putnam complex \cite{AP}. The
\v{C}ech cohomology of the hull of a primitive substitution tiling is
obtained as the direct limit of the substitution action on the
\v{C}ech cohomology of a finite CW-complex that approximates the
tiling space.  Let $A^{(m)}$ be the matrix of the
substitution action on the $m$-th cochain group of the approximant
complex. The zeta function is then given by \cite[Thm.~9.1]{AP} as
\begin{equation}
 \zeta(z) \, = \, \frac{\prod_{k\ \text{odd}}\det(1-zA^{(d-k)})}
                   {\prod_{k\ \text{even}}\det(1-zA^{(d-k)})}
   \, = \, \frac{\prod_{k\ \text{odd}}\prod_i(1-z\lambda_i^{(d-k)})}
                   {\prod_{k\ \text{even}}\prod_i(1-z\lambda_i^{(d-k)})}\ts ,
\label{zeta-ap}
\end{equation}
where the latter equality holds if all $A^{(m)}$ are diagonalizable,
with eigenvalues $\lambda_i^{(m)}$. Note that, instead of the action
on the cochain groups, one can also take the action on the cohomology
(with rational coefficients),
as the extra terms in (\ref{zeta-ap}) cancel between numerator and
denominator.

\begin{table}
  \caption{\label{tab:squiralcohom}
    Topology of the squiral hull and its various substitution factors. 
    In the left column, we indicate the set on which the 
    projection to $\mathbb{S}^{2}_{3}$ is not 1-to-1; 
    labels v or h mean that there are vertical or horizontal 
    sub-solenoids of type $\mathbb{S}^{1}_{3}$ where the 
    projection is 2-to-1, and an integer denotes the extra degeneracy 
    (beyond the 1d and 2d solenoids) of the projection to the origin
    of $\mathbb{S}^{2}_{3}$. In the second column, the cohomology 
    group $H^2$ is given. In all cases, $H^1=\ZZ[\frac{1}{3}]^2$ 
    and $H^0=\ZZ$. For the factors of $F_{\text{max}}$,
    examples for identifications of the symbols of $F_{\text{max}}$ are given in
    the last column. Especially for the smaller factors, there are usually many
    choices of identifications that yield equivalent factors.
  } 
\renewcommand{\arraystretch}{1.3}
\begin{tabular}{|l|l|l|}
\hline
multiplicity & $H^2$ & name / identifications \\
\hline
2-to-1 a.e. & $\ZZ[\frac19] \oplus  \ZZ[\frac13]^2 \oplus \ZZ^6$  & squiral \\
\hline
v,h,4     & $\ZZ[\frac19]\oplus\ZZ[\frac13]^2\oplus\ZZ^2\oplus\ZZ_2$ & 
$F_{\text{max}}$ \\
\hline
v,4       & $\ZZ[\frac19]\oplus\ZZ[\frac13]\oplus\ZZ^3$  & $f=g$ \\ 
v,3       & $\ZZ[\frac19]\oplus\ZZ[\frac13]\oplus\ZZ^2$  & 
$f=g$, $a=b$ or $c=d$ \\ 
v,2       & $\ZZ[\frac19]\oplus\ZZ[\frac13]\oplus\ZZ^1$  & 
$f=g$, $a=b$, $c=d$ \\ 
\hline
h,4       & $\ZZ[\frac19]\oplus\ZZ[\frac13]\oplus\ZZ^3$  & $e=g$ \\ 
h,3       & $\ZZ[\frac19]\oplus\ZZ[\frac13]\oplus\ZZ^2$  & 
$e=g$, $a=d$ or $b=c$ \\ 
h,2       & $\ZZ[\frac19]\oplus\ZZ[\frac13]\oplus\ZZ^1$  & 
$e=g$, $a=d$, $b=c$ \\ 
\hline
4         & $\ZZ[\frac19]\oplus\ZZ^4$                    & 
$e=f=g$ \\ 
3         & $\ZZ[\frac19]\oplus\ZZ^3$                    & 
$e=f=g$, $a=b$ \\ 
2         & $\ZZ[\frac19]\oplus\ZZ^2$                    & 
$e=f=g$, $a=b=c$ \\ 
1         & $\ZZ[\frac19]\oplus\ZZ^1$                    & 
$e=f=g$, $a=b=c=d$ \\ 
\hline
\end{tabular}
\end{table}

For the cohomology of the squiral tiling, we obtain
\begin{equation}
   H^2 \, =\, 
    \ZZ[\tfrac{1}{9}]\oplus\ZZ[\tfrac{1}{3}]^2\oplus\ZZ^6,\quad 
             H^1=\ZZ[\tfrac{1}{3}]^2, \quad  H^0=\ZZ\ts ,
\label{eq:cohomSq}
\end{equation} 
while the factor $F_{\text{max}}$ leads to
\begin{equation}
   H^2 \, =\, 
    \ZZ[\tfrac{1}{9}]\oplus\ZZ[\tfrac{1}{3}]^2\oplus\ZZ^2\oplus\ZZ_2,\quad 
             H^1=\ZZ[\tfrac{1}{3}]^2, \quad  H^0=\ZZ\ts .
\label{eq:cohomFmax}
\end{equation} 
In particular, the cohomology of $F_\text{max}$ contains
torsion. Taking into account the eigenvalues with which the
substitution acts on the cohomology, the zeta functions become
\begin{equation} 
   \zeta^{}_\text{sq}(z) \, = \, 
   \frac{(1-3z)^2}{(1-z)(1-9z)(1-3z)^2(1-z)^3(1+z)^3} \, =\,  
   \frac{1}{(1-z)(1-9z)(1-z^2)^3} \ts ,
\end{equation}
with fixed point numbers $a^{\text{(sq)}}_{m} = 9^m + 4 + 3\cdot (-1)^m$, 
and
\begin{equation}\label{f-zeta}
\begin{split} 
    \zeta^{}_{F_\text{max}}(z) \, & = \, 
    \frac{(1-3z)^2}{(1-z)(1-9z)(1-3z)^2(1-z)^2} 
    \, =\, \frac{1}{(1-z)^{3}(1-9z)} \\
     & = \,
    \frac{(1-3z)^2}{(1-z)(1-9z)} \cdot \left(\frac{1-z}{1-3z}\right)^2 
    \cdot \frac{1}{(1-z)^4} \ts ,
\end{split}
\end{equation}
with $a^{(F_{\text{max}})}_{m} = 9^m + 3$.  In the second line of
\eqref{f-zeta}, we have already expressed the zeta function of
$F_\text{max}$ as a product of the zeta function for
$\mathbb{S}^{2}_{3}$ (where $a^{(2)}_{m} = (3^m -1)^{2}$) the square
of the zeta function for $\mathbb{S}^{1}_{3}$ (with $a^{(1)}_{m} = 3^m
- 1$; compare sequence \texttt{A\ts 024023} of \cite{OEIS}) and the
zeta function for four extra fixed points. The corresponding zeta
functions follow from \cite{BLP}. Alternatively, since solenoids are
inverse limit spaces, the cohomology can also be computed with the
method from \cite{AP}.  The zeta function of $F_\text{max}$ clearly
reflects the structure of $F_\text{max}$ already determined
earlier. In contrast, the zeta function of the squiral cannot be
interpreted so easily.

There is a large number of further factors between $F_\text{max}$ and
the solenoid $\mathbb{S}^{2}_{3}$. We have systematically studied
those factors which can be obtained by identifying some of the seven
symbols of the substitution (\ref{eq:factorsubst}), while still giving
rise to a consistent substitution.  Such identifications of symbols
induce identifications of certain tilings projecting to the same point
on $\mathbb{S}^{2}_{3}$.

Under the identification of symbols $f=g$, pairs of tilings that
project to the same point of the vertical sub-solenoid are identified,
so that the splitting of preimages above that sub-solenoid is
closed. What remains is a split horizontal sub-solenoid, and four
extra fixed points.  The cohomology group $H^2$ loses the torsion
part, and one term $\ZZ[\frac13]$ is replaced by $\ZZ$. Closing the
splitting of the horizontal sub-solenoid with the identification $e=g$
is completely analogous.  Under the identification $e=f=g$, both 1d
sub-solenoid splittings are closed. Only four extra orbits above the
origin remain, and $H^2 = \ZZ[\frac{1}{9}] \oplus \ZZ^4$.  From here,
there are many ways to make further identifications. With each such
identification, the multiplicity of the projection to the origin of
$\mathbb{S}^{2}_{3}$ is reduced by 1, and so is the exponent of the
$\ZZ^k$ term in $H^2$. An example of a smallest factor above
$\mathbb{S}^{2}_{3}$ is obtained with the identifications $a=b=c=d$
and $e=f=g$, with $H^2=\ZZ[\frac19]\oplus\ZZ$.  Instead of closing the
horizontal split solenoid right after the vertical one, it is also
possible to begin with one or two of the identifications $a=b$ and
$c=d$. In such an operation, the exponent of the $\ZZ^k$ term of $H^2$
drops by 1. Examples of MLD classes of the different factors are
listed in Table~\ref{tab:squiralcohom}.

\section{The table tiling}\label{sec:table}

\begin{figure}[b]
\centerline{\includegraphics[width=0.4\textwidth]{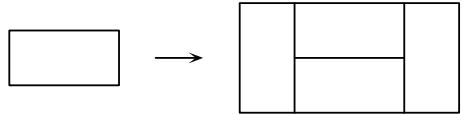}}
\caption{\label{tabinfl} Geometric inflation rule for the table tiling.}
\end{figure}

\begin{figure}
\centerline{\includegraphics[width=\textwidth]{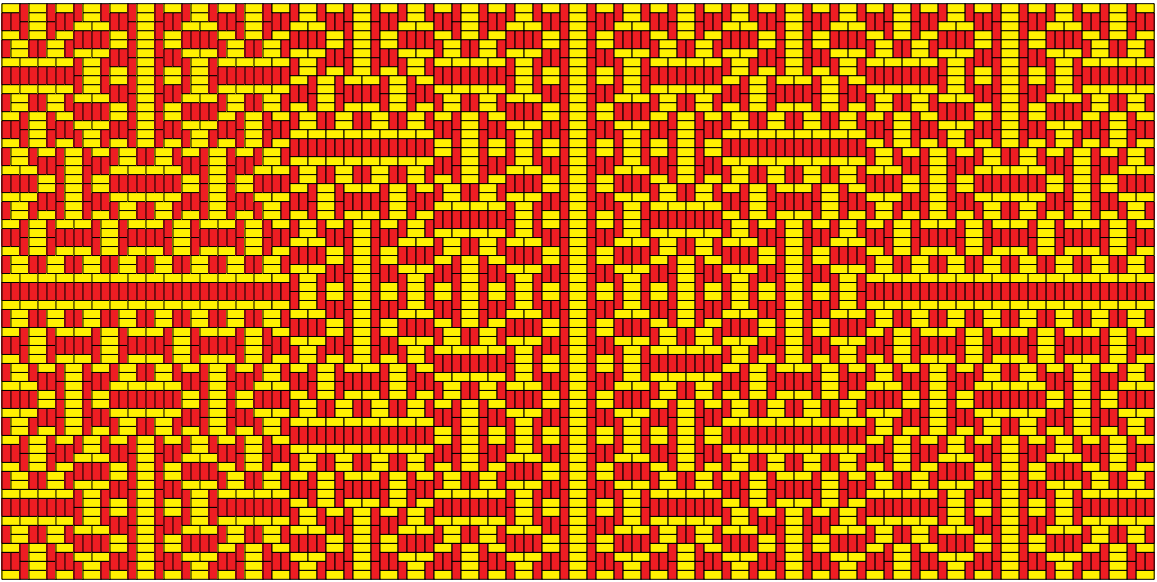}}
\caption{\label{tabpatch} A patch of the table tiling in the 
form of a level-$6$ supertile.}
\end{figure}

Infinite legal patterns generated by the block substitution
\begin{equation}\label{eq:blocktable}
   0\,\mapsto\,\left[\begin{array}{@{}c@{\;\;}c@{}}1&0\\[-1mm] 
     3&0\end{array}\right]\qquad
   1\,\mapsto\,\left[\begin{array}{@{}c@{\;\;}c@{}}0&2\\[-1mm] 
     1&1\end{array}\right]\qquad
   2\,\mapsto\,\left[\begin{array}{@{}c@{\;\;}c@{}}2&1\\[-1mm] 
     2&3\end{array}\right]\qquad
   3\,\mapsto\,\left[\begin{array}{@{}c@{\;\;}c@{}}3&3\\[-1mm] 
     0&2\end{array}\right]
\end{equation}
are MLD to the well-known table tilings with mixed spectrum
\cite{Rob99}. The geometric inflation rule is explained in
Figure~\ref{tabinfl}, and a larger patch is shown in
Figure~\ref{tabpatch}. Frettl\"oh \cite{Fre02} has independently shown
that the table tilings cannot be described as model sets. Indeed, from
the substitution (\ref{eq:blocktable}), one can see that the four
supertiles differ in every place, which is a generalisation of the
bijectivity property. As a consequence, table tilings that consist of
a single infinite order supertile occur in groups of four, which are
pairwise different at every position, and which all project to the
same point in the underlying solenoid $\mathbb{S}^2_2$. This means
that the projection to $\mathbb{S}^{2}_{2}$ is 4-to-1 almost
everywhere. In such a situation, the spectrum must indeed be mixed,
and the tiling space cannot be MLD to the hull of a model set.

What can be said about possible model set factors? In the case of
the Thue-Morse and squiral tilings, we have seen that there exist
maximal model set factors, onto which the projection has uniform
multiplicity. For the table tiling, such a factor with uniform
multiplicity cannot exist, which can be seen as follows.

For such a factor to exist, the multiplicity of the projection
to $\mathbb{S}^{2}_{2}$ must be divisible by four \emph{everywhere}.
However, there are 10 different pairs of tiles sharing a vertical
edge, and 10 pairs sharing a horizontal edge, so that for tilings
which consist of a pair of infinite order supertiles, sharing a
vertical or a horizontal line, the projection to $\mathbb{S}^{2}_{2}$ is 
10-to-1. As a consequence, any projection to a model set factor must 
identify all 10 tilings in such a group of 10, and thus cannot
have constant multiplicity. 

This does not exclude the existence of model set factors between the
table and the solenoid $\mathbb{S}^{2}_{2}$, however. Such a factor
can still have a projection to the solenoid with a non-trivial
multiplicity at the corners of infinite order supertiles.

The table tiling admits $24$ legal $2\!\times\! 2$-patterns, all of 
which are seeds for one of the $24$ fixed points under the square of 
the substitution. These patterns are the following,
\[
  \begin{split}
  \left[\begin{array}{@{}c@{\;\;}c@{}} 0 & 2 \\[-1mm] 
0 & 2 \end{array}\right] \;\;
  \left[\begin{array}{@{}c@{\;\;}c@{}} 0 & 2 \\[-1mm] 
1 & 0 \end{array}\right] \;\;
  \left[\begin{array}{@{}c@{\;\;}c@{}} 0 & 2 \\[-1mm]
 2 & 1 \end{array}\right] \;\;
  \left[\begin{array}{@{}c@{\;\;}c@{}} 1 & 0 \\[-1mm] 
3 & 1 \end{array}\right] \;\;
  \left[\begin{array}{@{}c@{\;\;}c@{}} 1 & 1 \\[-1mm]
 3 & 3 \end{array}\right] \;\;
  \left[\begin{array}{@{}c@{\;\;}c@{}} 1 & 3 \\[-1mm] 
3 & 0 \end{array}\right] \;\;
  \left[\begin{array}{@{}c@{\;\;}c@{}} 2 & 0 \\[-1mm]
 1 & 1 \end{array}\right] \;\;
  \left[\begin{array}{@{}c@{\;\;}c@{}} 2 & 1 \\[-1mm]
 1 & 3 \end{array}\right] \\[1mm]
  \left[\begin{array}{@{}c@{\;\;}c@{}} 2 & 3 \\[-1mm] 
0 & 2 \end{array}\right] \;\;
  \left[\begin{array}{@{}c@{\;\;}c@{}} 2 & 3 \\[-1mm] 
1 & 0 \end{array}\right] \;\;
  \left[\begin{array}{@{}c@{\;\;}c@{}} 2 & 3 \\[-1mm] 
2 & 1 \end{array}\right] \;\;
  \left[\begin{array}{@{}c@{\;\;}c@{}} 3 & 0 \\[-1mm] 
0 & 2 \end{array}\right] \;\;
  \left[\begin{array}{@{}c@{\;\;}c@{}} 3 & 0 \\[-1mm] 
1 & 0 \end{array}\right] \;\;
  \left[\begin{array}{@{}c@{\;\;}c@{}} 3 & 0 \\[-1mm] 
2 & 1 \end{array}\right] \;\;
  \left[\begin{array}{@{}c@{\;\;}c@{}} 3 & 1 \\[-1mm] 
2 & 3 \end{array}\right] \;\;
  \left[\begin{array}{@{}c@{\;\;}c@{}} 3 & 3 \\[-1mm]
 2 & 0 \end{array}\right] \\[1mm]
  \left[\begin{array}{@{}c@{\;\;}c@{}} 0 & 2 \\[-1mm] 
2 & 0 \end{array}\right] \;\;
  \left[\begin{array}{@{}c@{\;\;}c@{}} 1 & 3 \\[-1mm] 
3 & 1 \end{array}\right] \;\;
  \left[\begin{array}{@{}c@{\;\;}c@{}} 2 & 0 \\[-1mm] 
0 & 2 \end{array}\right] \;\;
  \left[\begin{array}{@{}c@{\;\;}c@{}} 3 & 1 \\[-1mm] 
1 & 3 \end{array}\right] \;\;
  \left[\begin{array}{@{}c@{\;\;}c@{}} 0 & 2 \\[-1mm] 
1 & 1 \end{array}\right] \;\;
  \left[\begin{array}{@{}c@{\;\;}c@{}} 1 & 0 \\[-1mm] 
3 & 0 \end{array}\right] \;\;
  \left[\begin{array}{@{}c@{\;\;}c@{}} 2 & 1 \\[-1mm] 
2 & 3 \end{array}\right] \;\;
  \left[\begin{array}{@{}c@{\;\;}c@{}} 3 & 3 \\[-1mm] 
0 & 2 \end{array}\right]
   \end{split}
\]
We have found two sliding block maps which induce a primitive 
substitution on the resulting factor. In the first map, the
$16$ patches of the first two rows map to $0$ and the remaining 
eight patches from the third row to $1$. This factor map induces 
a block substitution
\begin{equation}
    \ell\,\mapsto\,\left[\begin{array}{@{}c@{\;\;}c@{}}0&\ell\\[-1mm] 
                     1&0\end{array}\right]\qquad
    \label{eq:tablefac1subst}
\end{equation}
on the alphabet $\{0,1\}$, with $\ell\in\{0,1\}$. The new substitution
rule is primitive and admits two fixed points (with legal seed
$\left[\begin{smallmatrix} 0 & 1 \\ 1 & 0 \end{smallmatrix}\right]$ or
$\left[\begin{smallmatrix} 0 & 1 \\ 0 & 0 \end{smallmatrix}\right]$
and reference point in its center, so that the fixed point covers
$\ZZ^2$).  The substitution obviously has a Toeplitz structure, so
that each legal tiling it generates forms a model set, with pure point
spectrum; compare \cite{LM,FS}.  The two fixed points project to the
same point in $\mathbb{S}^{2}_{2}$, so that there are two translation
orbits of tilings which project to a single translation orbit in
$\mathbb{S}^{2}_{2}$.

Another (larger) factor is obtained if the first $16$ patches are
mapped to $0$, the first $4$ patches in the third row to $1$, and the
remaining $4$ patches to $2$. The induced substitution on the alphabet
$\{0,1,2\}$ is
\begin{equation}
    0\,\mapsto\,\left[\begin{array}{@{}c@{\;\;}c@{}}2&1\\[-1mm]
 1&2\end{array}\right]\, ,\qquad
    1\,\mapsto\,\left[\begin{array}{@{}c@{\;\;}c@{}}2&0\\[-1mm]
 1&2\end{array}\right]\, ,\qquad
    2\,\mapsto\,\left[\begin{array}{@{}c@{\;\;}c@{}}2&2\\[-1mm] 
1&2\end{array}\right]\, .
    \label{eq:tablefac2subst}
\end{equation}
Again, the substitution is primitive and has a Toeplitz structure,
implying that also this factor is a model set, with pure point
spectrum \cite{LM,FS}.  There are three fixed points under the square
of the substitution, with seeds $\left[\begin{smallmatrix} 1 & 2 \\
    \ell & 1 \end{smallmatrix}\right]$ ($\ell=1,2,3$), which project
to the same point on $\mathbb{S}^{2}_{2}$, so that there are
three translation orbits of tilings which project to a single
translation orbit on the solenoid. Everywhere else, the projection is
1-to-1.

The interpretation of the zeta function turns out to be more difficult
if the projection to the solenoid fails to be 1-to-1 almost
everywhere. For the table tiling, we obtain as \v{C}ech cohomology
\begin{equation}
   H^2 \, =\, 
    \ZZ[\tfrac{1}{4}]\oplus\ZZ[\tfrac{1}{2}]^4\oplus\ZZ^3\oplus\ZZ_2,\quad 
             H^1=\ZZ[\tfrac{1}{2}]^2, \quad  H^0=\ZZ\ts ,
\label{eq:cohomTable}
\end{equation} 
while the dynamical zeta fuction reads
\begin{equation}\label{table-zeta}
  \zeta^{}_\text{table}(z)\, = 
    \,\frac{(1-2z)^2}{(1-4z)(1-4z^2)^2(1-z)^2(1-z^2)} \ts .
\end{equation}
The corresponding fixed point counts are $a^{}_{m} = 4^n + 3 + (-1)^n
(1 + 2^{n+1})$.  It is natural to compare this to the zeta function of
$\mathbb{S}_2^2$,
\begin{equation}\label{table-s22}
  \zeta^{}_{\mathbb{S}_2^2}(z) \,  = \, \frac{(1-2z)^2}{(1-4z)(1-z)} \ts , 
\end{equation}
which follows once again from \cite{BLP} and codes the fixed
point counts $a^{(2)}_{m} = (2^n - 1)^2$; compare sequences
\texttt{A\ts 000225} and \texttt{A\ts 060867} of \cite{OEIS}.
At this stage, however, it is far from obvious how the fixed
point counts of the table shall be mapped to those of the solenoid,
under a map which is 4-to-1 almost everywhere.

\section*{Acknowledgement}
It is a pleasure to thank Jean-Paul Allouche for many fruitful
discussions. Moreover, we are grateful to Natalie Frank, Aernout van
Enter and Daniel Lenz for helpful discussions on the example from
Section~\ref{sec:fun}, which took place during a meeting at BIRS in
November 2011.  This work was supported by the German Research Council
(DFG), within the CRC 701.


\begin{thebibliography}{99}

\bibitem{Apal}
Allouche J-P,
Schr\"{o}dinger operators with Rudin-Shapiro potentials
are not palindromic,
\textit{J.\ Math.\ Phys.} \textbf{38} (1997) 1843--1848.

\bibitem{ABCD}
Allouche J-P, Baake M, Cassaigne J and Damanik D,
Palindrome complexity, 
\textit{Theor.\ Comput.\ Science} \textbf{292} (2003) 9--31;
\texttt{arXiv:math.CO/0106121}.

\bibitem{AB}
Allouche J-P, Bacher R,
Toeplitz sequences, paperfolding, Hanoi towers and progression-free
sequences of integers,
\textit{Ens.\ Math.} \textbf{38} (1992) 315--327.

\bibitem{AMF}
Allouche J-P and Mend\`{e}s France M,
Automatic sequences,
in:\ \textit{Beyond Quasicrystals},
eds Axel F and Gratias D,
Springer, Berlin (2000) pp.~293--367.

\bibitem{AS93}
Allouche J-P and Shallit J,
Complexit\'{e} des suites de Rudin-Shapiro g\'{e}n\'{e}ralis\'{e}es,
\textit{J.\  Th\'{e}orie Nombres de Bordeaux}
\textbf{5} (1993) 283--302.

\bibitem{AS}
Allouche J-P and Shallit J,
\textit{Automatic Sequences:\ Theory,\ Applications,\ Generalizations},
Cambridge University Press, Cambridge (2003).

\bibitem{AP}
Anderson J E and Putnam I F,
Topological invariants for substitution tilings and their 
associated \break $C^{\ast}$-algebras,
\textit{Ergodic Th.\ \& Dynam.\ Syst.} \textbf{18} (1998) 509--537.

\bibitem{Bpal}
Baake M,
A note on palindromicity,
\textit{Lett.\ Math.\ Phys.} \textbf{49} (1999) 217--227;
\texttt{arXiv:math-ph/9907011}.

\bibitem{B02}
Baake M,
A guide to mathematical quasicrystals,
in:\ \textit{Quasicrystals -- An Introduction to Structure,
Physical Properties and Applications},
eds J-B Suck, M Schreiber and  P H\"{a}ussler,
Springer, Berlin (2002) pp.~17--48;
\texttt{arXiv:math-ph/9901014}.

\bibitem{BGG}
Baake M, G\"{a}hler F and Grimm U,
Spectral and topological properties of a family of 
generalised Thue-Morse sequences,
\textit{J.\ Math.\ Phys.} \textbf{53} (2012) 032701;
\texttt{arXiv:1201.1423}. 

\bibitem{BGrev}
Baake M and Grimm U,
Kinematic diffraction from a mathematical viewpoint,
\textit{Z.\ Kristallogr.} \textbf{226} (2011) 711--725;
\texttt{arXiv:1105.0095}.

\bibitem{BEL} 
Baake M, Lenz D and van Enter A C D, 
Diffraction versus dynamical spectra for uniquely ergodic systems with
finite local complexity, in preparation.

\bibitem{BE}
Baake M and van Enter A C D,
Close-packed dimers on the line:\ Diffraction
versus dynamical spectrum,
\textit{J.\ Stat.\ Phys.} \textbf{143} (2011) 88--101;
\texttt{arXiv:1011.1628}.

\bibitem{BG12}
Baake M and Grimm U,
Squirals and beyond:\ Substitution tilings with singular continuous
spectrum, \textit{Ergodic Th.\ \& Dynam.\ Syst.}, to appear;
\texttt{arXiv:1205.1384}.

\bibitem{TAO}
Baake M and Grimm U,
\textit{Theory of Aperiodic Order:\ A Mathematical Invitation},
Cambridge University Press, in preparation.

\bibitem{BLP}
Baake M, Lau E and Paskunas V,
A note on the dynamical zeta function of general
toral endomorphisms,
\textit{Monatsh.\ Math.} \textbf{161} (2010) 33--42;
\texttt{arXiv:0810.1855}.

\bibitem{BL}
Baake M and Lenz D,
Dynamical systems on translation bounded measures:\
Pure point dynamical and diffraction spectra, 
\textit{Ergodic Th.\ \& Dynam.\ Syst.} \textbf{24} (2004)
1867--1893; \texttt{arXiv:math.DS/0302231}.

\bibitem{BLM}
Baake M, Lenz D and Moody R V,
Characterization of model sets by dynamical systems, 
\textit{Ergodic Th.\ \& Dynam.\ Syst.} \textbf{27} (2007) 341--382;
\texttt{arXiv:math.DS/0511648}.

\bibitem{BMBook}
Baake M and Moody R V (eds.),
\textit{Directions in Mathematical Quasicrystals},
CRM Monograph Series vol.\ 13, 
AMS, Providence, RI (2000).

\bibitem{BM}
Baake M and Moody R V,
Weighted Dirac combs with pure point diffraction,
\textit{J.\ reine angew.\ Math.\ (Crelle)} 
\textbf{573} (2004) 61--94;
\texttt{arXiv:math.MG/0203030}.

\bibitem{BSJ}
Baake M, Schlottmann M and Jarvis P D, 
Quasiperiodic patterns with tenfold symmetry and equivalence with 
respect to local derivability,  
\textit{J.\ Phys.~A:\ Math.\ Gen.} \textbf{24} (1991) 4637--4654.

\bibitem{Cow}
Cowley J M,
\textit{Diffraction Physics},
3rd ed., North-Holland, Amsterdam (1995).

\bibitem{D}
Dekking F M,
The spectrum of dynamical systems arising from substitutions 
of constant length,
\textit{Z.\ Wahrscheinlichkeitsth.\ verw.\ Geb.} \textbf{41} 
(1978) 221--239.

\bibitem{EM}
van Enter A C D and Mi\c{e}kisz J,
How should one define a weak crystal?
\textit{J.\ Stat.\ Phys.} \textbf{66} (1992) 1147--1153.

\bibitem{Nat1}
Frank N P,
Multi-dimensional constant-length substitution sequences,
\textit{Topol.\ Appl.} \textbf{152} (2005) 44--69.

\bibitem{Nat2}
Frank N P,
Spectral theory of bijective substitution sequences,
\textit{MFO Reports} \textbf{6} (2009) 752--756.

\bibitem{Fre02}
Frettl\"{o}h D.\ (2002).
\textit{Nichtperiodische Pflasterungen mit ganzzahligem Inflationsfaktor},
PhD thesis (Univ.\ Dortmund).

\bibitem{FS}
Frettl{\"o}h D and Sing B,
Computing modular coincidences for substitution tilings and point sets,
\textit{Discr.\ Comput.\ Geom.} \textbf{37} (2007) 381--407;
\texttt{arXiv:math.MG/0601067}.

\bibitem{GS}
Gr\"{u}nbaum B and Shephard G C,
\textit{Tilings and Patterns}, Freeman, New York (1987).

\bibitem{Hof}
Hof A,
On diffraction by aperiodic structures,
\textit{Commun.\ Math.\ Phys.} \textbf{169} (1995) 25--43. 
 
\bibitem{Kaku} 
Kakutani S,
Strictly ergodic symbolic dynamical systems,
\textit{Proc.\ 6th Berkeley Symposium on Math.\ Statistics
and Probability} eds L M LeCam, J Neyman and 
E L Scott, Univ.\ of California Press, Berkeley (1972),
pp.\ 319--326.

\bibitem{Kea}
Keane M, 
Generalized Morse sequences,
\textit{Z.\ Wahrscheinlichkeitsth.\  verw.\ Geb.}  
\textbf{10} (1968) 335--353.

\bibitem{Kit}
Kitchens B P,
\textit{Symbolic Dynamics},
Springer, Berlin (1998).

\bibitem{LM}
Lee J-Y and Moody R V,
Lattice substitution systems and model sets, 
\textit{Discr.\ Comput.\ Geom.} \textbf{25} (2001) 173--201;
\texttt{arXiv:math.MG/0002019}.

\bibitem{LMS}
Lee J-Y, Moody R V and Solomyak B,
Pure point dynamical and diffraction spectra,
\textit{Ann.\ H.\ Poincar\'{e}} \textbf{3} (2002) 1003--1018;
\texttt{arXiv:0910.4809}. 

\bibitem{LMBook}
Lind D A and Marcus B (1995),
\textit{An Introduction to Symbolic Dynamics and Coding},
Cambridge University Press, Cambridge (1995).

\bibitem{Moody00}
Moody R V, 
Model sets:\ A Survey,
in \textit{From Quasicrystals to More Complex Systems},
Axel F, D\'enoyer F and Gazeau J P  (eds.),
EDP Sciences, Les Ulis, and Springer, Berlin (2000), 
pp.\ 145--166;
\texttt{arXiv:math.MG/0002020}.

\bibitem{OEIS}
\textit{The On-Line Encyclopedia of Integer Sequences}$^{\text{TM}}$,
\texttt{http://oeis.org/}.

\bibitem{PF}
Pytheas Fogg N,
\textit{Substitutions in Dynamics, Arithmetics
and Combinatorics}, LNM 1794,
Springer, Berlin (2002).

\bibitem{Q}
Queff\'{e}lec M,
\textit{Substitution Dynamical Systems -- Spectral Analysis},
LNM 1294, 2nd ed., Springer, Berlin (2010).

\bibitem{Rob99}
Robinson E A,
On the table and the chair,
\textit{Indag.\ Math.} \textbf{10} (1999) 581--599.

\bibitem{Rob04}
Robinson E A,
Symbolic dynamics and tilings of $\RR^{d}$,
\textit{Proc.\ Sympos.\ Appl.\ Math.} \textbf{60} (2004) 81--119.

\bibitem{Ruelle}
Ruelle D,
\textit{Dynamical Zeta Functions for Piecewise Monotone 
Maps of the Interval}, CRM Monograph Series, vol.\ 4, 
AMS, Providence, RI (1994).

\bibitem{Sadun}
Sadun L,
\textit{Topology of Tiling Spaces}, 
AMS, Providence, RI (2008).

\bibitem{Schmidt}
Schmidt K,
\textit{Dynamical Systems of Algebraic Origin},
Birkh\"{a}user, Basel (1995).

\bibitem{Danny}
Shechtman D, Blech I, Gratias D and Cahn J W,
Metallic phase with long-range orientational order and no
translational symmetry,
\textit{Phys.\ Rev.\ Lett.} \textbf{53} (1984) 1951--1953.

\bibitem{M}
Schlottmann M,
Generalised model sets and dynamical systems,
in:\ \cite{BMBook}, pp.\ 143--159.

\bibitem{Withers}
Withers R L,
Disorder, structured diffuse scattering and the
transmission electron microscope,
\textit{Z.\ Krist.} \textbf{220} (2005) 1027--1034.

\end{thebibliography}
\end{document}